\documentclass[%12pt
]{amsart}

\usepackage{amssymb,amscd}
\usepackage{enumerate,graphicx}

\theoremstyle{plain}
\newtheorem{theorem}{Theorem}
\newtheorem{lemma}[theorem]{Lemma}
\newtheorem{proposition}[theorem]{Proposition}
\newtheorem{corollary}[theorem]{Corollary}

\theoremstyle{definition}
\newtheorem{definition}[theorem]{Definition}

\theoremstyle{remark}
\newtheorem*{remark}{Remark}
\newtheorem*{remarks}{Remarks}

\newtheorem*{example}{{\bf Example}}

\newcommand{\Aa}{\hbox{\ensuremath{\mathcal{A}}}}
\newcommand{\HH}{\hbox{\ensuremath{\mathcal{H}}}}
\newcommand{\Ss}{\hbox{\ensuremath{\mathcal{S}}}}
\newcommand{\KK}{\hbox{\ensuremath{\mathcal{K}}}}
\newcommand{\QQ}{Q}
\newcommand{\BQ}{\hbox{\ensuremath{\mathcal{B}_{Q}}}}
\newcommand{\KQ}{\hbox{\ensuremath{\mathcal{K}_{Q}}}}

\newcommand{\KS}{\hbox{\ensuremath{\mathcal{K}_{S}}}}
\newcommand{\KV}{\hbox{\ensuremath{\mathcal{K}_{V}}}}

\newcommand{\Z}{\hbox{\ensuremath{\mathbb{Z}}}}

\newcommand{\R}{\hbox{\ensuremath{\mathbb{R}}}}
\newcommand{\C}{\hbox{\ensuremath{\mathbb{C}}}}
\newcommand{\supess}{\hbox{\ensuremath{\text{ess}\sup}}}
\newcommand{\card}{\hbox{\ensuremath{\#}}}
\newcommand{\comment}[1]{}

\newcommand{\mnorm}[1]{%
  \left\vert\kern-0.9pt\left\vert\kern-0.9pt\left\vert #1
    \right\vert\kern-0.9pt\right\vert\kern-0.9pt\right\vert}
\DeclareMathOperator{\supp}{Supp}

\numberwithin{equation}{section}
\numberwithin{theorem}{section}

\begin{document}

%%%%%%%%%%%%%%%%%%%%%%%% PREAMBULO

%\SetTexturesEPSFSpecial
%\HideDisplacementBoxes

%%%%%%%%%%%%%A continuacion va el titulo

\title[Wavelet frames] {Wavelets on Irregular Grids with Arbitrary
Dilation Matrices, and
Frame Atoms for $L^2(\R^d)$}

% Author Info

\author[A.Aldroubi]{Akram Aldroubi} \address{Department of Mathematics\\,
Vanderbilt University\\ 1326 Stevenson Center\\Nashville, TN 37240}
\email[Akram Aldroubi]{aldroubi@math.vanderbilt.edu}

\author[C.Cabrelli]{Carlos~Cabrelli}
\address{Departamento de
Matem\'atica \\ Facultad de Ciencias Exactas y Naturales\\ Universidad
de Buenos Aires\\ Ciudad Universitaria, Pabell\'on I\\ 1428 Capital
Federal\\ ARGENTINA\\ and CONICET, Argentina}
\email[Carlos~Cabrelli]{cabrelli@dm.uba.ar} \thanks{The research of
Akram Aldroubi is supported in part by NSF grant DMS-0103104, and by
DMS-0139740. The research of
       Carlos Cabrelli and Ursula Molter is partially supported by
Grants: PICT 03134, and CONICET, PIP456/98}

\author[U.Molter]{Ursula~M.~Molter}
\address{Departamento de
Matem\'atica \\ Facultad de Ciencias Exactas y Naturales\\ Universidad
de Buenos Aires\\ Ciudad Universitaria, Pabell\'on I\\ 1428 Capital
Federal\\ ARGENTINA\\ and CONICET, Argentina}
\email[Ursula~M.~Molter]{umolter@dm.uba.ar}

%%%%%%%%%%%%%%%%%%% Aqui van los keywords y la subject class


\date{\today}

%%%%%%%%%%%%% Abstract

\begin{abstract}
%texto del abstract aqui
In this article, we develop a general method for con\-struct\-ing
wa\-ve\-lets $\{|\det A_j|^{1/2}\psi(A_jx-x_{j,k}): \,
j\in J, k \in K\}$  on irregular lattices of the form
$X=\{x_{j,k} \in \R^d: \; j \in J, k\in K\}$, and with an arbitrary countable
family of invertible $d\times d$
matrices $\{A_j \in GL_d(\R): \; j \in J\}$ that do not necessarily
have a group
structure. This wavelet
construction is a particular case of  general atomic frame
decompositions of $L^2(\R^d)$ developed in this article, that allow
other time frequency decompositions such as non-harmonic Gabor frames with
non-uniform covering of the Euclidean space $\R^d$. Possible
applications include image and video compression, speech coding,
image and digital data transmission,
image analysis, estimations and detection, and seismology.

\end{abstract}

\maketitle

%%%%%%%%%%%%%%%%%%%%%%%%%% FIN PREAMBULO

%%%%%%%%%%%%%% A partir de aqui va el paper

\section{Introduction}
\label {I}
Recently there has been a considerable interest in trying to obtain
atomic decompositions
of the space $ L^2(\R^d)$. These decompositions are usually obtained in terms
of frames generated by a family of functions translated on a regular
grid, and dilated by powers
of a dilation matrix. The uniformity of the grid and the structure of
the dilations can be exploited to
obtain very sharp results. For irregular grids and unstructured
dilations or if dilations are replaced by other
transformations the situation is more complex and requires different
techniques.
One method is to use the regular case and try to obtain perturbations
of the grid that preserve the frame structure.
Another possibility is to obtain
irregular samples of the continuous transform, that have the
required properties.

In this article we study frame decompositions of the space $L^2(\R^d)$
using translations of a family of functions on irregular grids, and
arbitrary dilations, and we even
replace dilations by other transformations.

Our approach is different and very general, allowing quite general
constructions. We prove the existence of smooth time-frequency frame atoms in
several variables. The setting  includes as particular cases, wavelet
frames on irregular lattices and
with a set of dilations or transformations that do not have a group
structure. Another
particular case are non-harmonic Gabor frames with non-uniform
covering of the Euclidean space. It also
leads to new constructions of  wavelet and Gabor frames with regular
lattice translates. One of the nice
features of the proposed method is that it unifies  different atomic
decompositions.

For the case of regular lattices Guido Weiss and his group \cite{HLW02,%
HLWW02,Lab02} developed
a very fundamental program to characterize a large class of
decompositions of
$L^2(\R^d)$ through certain equations that the generators must satisfy.
This is
an important attempt to unify Gabor and wavelets decompositions.
Other fundamental construction of MRA wavelet
frames on regular lattices can also be found in \cite {CHS}, \cite
{CHSS03}, \cite {CS00}.
Our methods can be used to produce a substantial part of these systems.

A set $Q \subset \R^d$ is a wavelet set if the
inverse Fourier transform of the characteristic
function of the set is a wavelet. Wavelet sets,
frame wavelet sets and methods for constructing
such sets have been studied recently \cite
{BMM99}, \cite{BL99}, \cite{BL01}, \cite{BS03},
\cite{DLS97}, \cite{DLS98}, \cite {HL00}
    \cite{Ola03}, \cite{OS03}. Our methods give constructions of wavelet
sets with translations on
irregular grids.

Let $J$ and $K$ be countable index sets. We consider families of functions
$\{g_j\}_{j \in J} \subset L^2(\R^d)$ and
discrete sets $ X=\{x_{j,k}: j \in J, k \in K\} \subseteq \R^d$
such that the
collection
$ \left\{g_j(x-x_{j,k}):  j \in J, k \in K\right\}$
form a frame for $L^2(\R^d)$.
The wavelet case is obtained when $g_j= |\det(A)|^{j/2}  g \circ A^j$ with A an
expansive matrix and $g$ a fixed atom.
We want to stress here  that our constructions are much more general, allowing
for example a different invertible (not necessarily expansive) matrix $A_j$
for each $j \in J.$
For the case of orthogonal wavelets, Yang Wang  \cite{Wan02}
has recently considered
wavelet
sets associated with arbitrary families of invertible matrices
and irregular sets of translates.
He gave conditions for the existence of such wavelet sets
and related them to spectral pairs.

Irregular wavelet and Gabor frames also have been studied as perturbations of
uniform (lattice translate) frames and also as sampling of the continuous
wavelet/ Gabor transform. See \cite{ Bal97}, \cite {BCHL03}, \cite{ Chr96},
    \cite{ Chr97}, \cite{CFZ01},
    \cite {CH97}, \cite{ CDH99}, \cite{ FZ95}, \cite{ FG89},
    \cite{ Gro91}, \cite{ Gro93},\cite {HK03}, \cite{ OS92}, \cite{ RS95},
    \cite{ SZ00}, \cite{ SZ01}, \cite{ SZ02a}, \cite{ SZ03},
    \cite{ SZ03a}.

The approach in this article can be considered in the spirit
of the classic construction in 1 dimension of smooth regular tight frames done
by  Daubechies, Grossmann and Meyer in  \cite{DGM86}. They found, for
the case of uniform lattices,
general conditions on a compactly supported smooth  function $h$,
in order that it
generates a tight Gabor frame of $L^2(\R)$. In the affine case they
found necessary and sufficient
conditions for a band limited function in order that it generates a
smooth wavelet frame.  See also \cite{HW89}.

There were other related attempts to obtain atomic decompositions of
functional spaces using very general systems. See for example
\cite{FG85},
\cite{Fei87} in the context of locally compact groups.

This paper is organized as follows: Section \ref {notation}
introduces the notation and some
preliminaries. Section \ref {WC} presents a Theorem on wavelet
construction on arbitrary,
sufficiently dense, but otherwise irregular grids and with
arbitrary dilation or even invertible
transformation matrices. Specific constructions of such wavelets are
obtained in Section \ref {EW},
first in the 1-D case and then in the multidimensional case. A
general theory of frame atomic
decomposition of $L^2(\R^d)$ is obtained in Section \ref {gen-res}.
Using the concept of outer
frame, reconstruction formulas for these atomic decompositions are
obtained in Section \ref {RF}.

\section{Notation}
\label{notation}

Throughout the paper $J$ and $K$ will denote countable index sets,
and $e_x$ will stand for the function
$e_{x}(\xi) = e^{-i2\pi x\cdot\xi}.$ We will use $\mu (E)$ to denote
the Lebesgue
measure of a measurable set $E$.

A set $\mathcal H := \{h_j\}_{j\in J}$ of measurable functions on $
\R^d$ is called
a {\em Riesz partition of unity} ({\bf RPU}), if
there exist constants $0<p\leq P < +\infty$ such that
\begin{equation}\label{riesz-partition}
p \leq \sum_{j\in J} | h_j(x) |^2 \leq P \quad \text{a.e.}\quad x \in \R^d.
\end{equation}

Let $\Ss = \{S_j\}_{j \in J}$ be a family of measurable subsets  of
$\R^d$.
A {\em Riesz partition of unity {\bf associated to $\Ss$}}, is
a set $\mathcal H := \{h_j\}_{j\in J}$ of measurable functions, such that
\begin{enumerate}
\item $\supp h_j \subseteq S_j$
\item There exist constants $0<p\leq P < +\infty$ such that
\begin{equation}
p \leq \sum_{j\in J} | h_j(x)|^2 \leq P \quad \text{a.e.}\quad x \in
\cup_j S_j.
\end{equation}
\end{enumerate}

\begin{remarks}
\
\begin{itemize}
\item If $\mathcal H = \{h_j\}$ is a RPU, then
$\overline{\mathcal H} = \{\overline h_j\}$ is also a RPU.
\item If $p=P=1$, we will say that $\mathcal H = \{h_j\}$
is a {\em regular partition of unity}.
\item
If the sets in $\Ss$ are essentially disjoint
(i.e. $\mu(S_i\cap S_j) = 0, \forall i \not= j$),
the family $\{h_j = \chi_{S_j}\}$ will yield a regular partition of unity
associated to $\Ss$.
\item Every RPU $\mathcal H = \{h_j\}$ can be normalized
to obtain a regular partition of unity by considering
$$
\tilde{h}_j = \frac{h_j}{(\sum_j |h_j|^2)^{1/2}}.
$$
\item Given a family $\Ss = \{S_j\}_{j\in J}$ of measurable sets on $\R^d$,
define
\begin{equation} \label{covering-index}
\rho_S(x) = \card(\{j \in J: x \in S_j\}) =
\sum_{j \in J} \chi_{S_j}(x),
\end{equation}
where $\card(B)$ is the cardinal of the set $B$.
The value $\rho_S = \|\rho_S\|_{\infty}$
is called the {\em covering index} of $\Ss$.
\end{itemize}
\end{remarks}

We now recall the definition of frame for a given close subspace $F$ of
$L^2(\R^d)$.

\begin{definition}
A set of functions $\{g_j\}_{j\in J}$ is a {\em frame} for $F$  if
$g_j \in F$  and there exist constants
$0<m, M < +\infty$, such that
\begin{equation}\label{frame}
m\|f\|^2 \leq \sum_{j \in J} |<f,g_j>|^2 \leq
M\|f\|^2, \quad \forall \quad f\in F.
\end{equation}
\end{definition}

For a measurable set $Q \subseteq \R^d$ we will denote by $\KQ$  the
functions that have
support in $\overline \QQ$, and  by $\BQ$ the functions whose Fourier transform
has support in $\overline \QQ$ , i.e.
\begin{align} \label{wienerq}
\KQ & :=\{f \in L^2(\R^d) :
\supp (f) \subset {\overline{\QQ}}\} \\
\BQ & :=\{f \in L^2(\R^d) :
\supp (\hat f) \subset {\overline {\QQ}}\}.
\label{wienereq-1}
\end{align}

We will need also the following definition.
\begin{definition}
Let $H$ be a Hilbert space.
A collection of vectors $\{g_j \in H\}_{j\in J}$ is an {\em outer frame}
for a closed subspace $F$ of $H$,  if $\{P_F(g_j)\}_{j \in J}$ is
a frame for $F$, where $P_F$ is the orthogonal projection onto
$F$, or equivalently,
there exist constants
$0<m, M < +\infty$, such that
\begin{equation}\label{outer-frame}
m\|f\|^2 \leq \sum_{j \in J} |<f,g_j>|^2 \leq
M\|f\|^2, \quad \forall \quad f\in F.
\end{equation}
\end{definition}

Related definitions to the concept of outer frames appear in
\cite{LO98}, \cite{FW01}.

\begin{remark} Throughout the paper, we will use the following
immediate and very useful fact about frames:\\
If $\{g_j\}_{j \in J}$ is a frame for $\KQ$, and $V \subset \QQ$, then
$\{g_j\}_{j \in J}$ is an outer frame for $\mathcal{K}_V$.
\end{remark}

\section{Wavelets on arbitrary irregular grids and with arbitrary
dilation matrices and other transformations}
\label {WC}
Our first results concerns the construction of
wavelets $\{|\det A_j|^{1/2}\psi(A_jx-x_{j,k}):\, j\in J, k\in K\}$
with translates
on the arbitrary irregular grid
$X=\{x_{j,k} \in \R^d: \; k\in K, j \in J\}$ and with an arbitrary countable
family of invertible $d\times d$
matrices $\{A_j \in GL_d(\R): \; j \in J\}$.
\begin {theorem} [Wavelets]
\label {wavelet}
Let $\QQ\subset \R^d$ be a set of finite measure, $h$ a
function in $L^2(\R^d)$ and $\Aa = \{A_j \in
GL_d(\R):\, j \in J\}$ a family of invertible matrices.

For each $j \in J$ set
$
B_j = (A_j^T)^{-1},\ S_j = B_j^{-1} \QQ = A_j^T \QQ,\ h_j = h(B_j\cdot)\
$
and let
$\Ss = \{S_j, j\in J\}$.

Assume that $\Ss$ is a covering of $\R^d$,
$\HH$ is a RPU with bounds $p$ and $P$ and that $\supp(h) 
\subset \QQ$.

Consider $X = \{x_{j,k} \in \R^d : j \in J, k \in K\}$ such that
for each $j \in J$, the set  $\{e_{x_{j,k}}\chi_{\QQ}: k \in K\}$
forms a frame for
$\mathcal{K}_{\QQ}$ with  lower and upper frame bounds $m_j$ and
$M_j$ respectively.
If $m := \inf_{j} m_j>0$ and
$M := \sup_j M_{j} < +\infty$, then the collection
\begin{equation*}
\{|\det A_j|^{1/2}\psi(A_jx-x_{j,k}):\, j\in J, k\in K\}
\end{equation*}
is a wavelet frame of $L^2(\R^d)$ with bounds $mp$ and $MP$,  generated by a single function $\psi$,
where $\psi$ is the inverse Fourier transform of $h$.
\end{theorem}

\begin{proof}
Since for each $j \in J$ we have that $\{e_{x_{j,k}} \chi_{\QQ} : k \in K \} $
  forms a frame for $\mathcal{K}_{\QQ}$ with
lower and upper frame bounds $m_j$ and $M_j$ respectively,
an application of Part~\ref{coro-2} of Corollary \ref{transl} for the
matrix $B_j^{-1}*$ shows that
$\{|B_j|^{1/2} e_{x_{j,k}} (B_j\omega) \chi_{\QQ} (B_j\omega) : k \in K \} $
forms a frame of ${\mathcal K}_{S_j}$ with the same bounds.
 From the definition of $S_j$,
$\{(\mu(S_j))^{- (1/2)} e_{(A_{j}^{-1}x_{j,k})}(\omega) 
\chi_{S_j}(\omega) : k \in K \}$
is then a frame for ${\mathcal K}_{S_j}$ with frame bounds 
$m_j\mu(\QQ)^{-1}$ and
$M_j\mu(\QQ)^{-1}$.

On the other side, if $\supp(h) \subset\QQ$ then $\HH$ is associated 
to $S$. So, we can apply Proposition \ref{genexpdil} to $S,\HH$ and
  $ \{\mu(S_j)^{-(1/2)} e_{(A_{j}^{-1}x_{j,k})} \chi_{S_j}: k \in K 
\}, j \in J$,
to conclude that,
$\{ \mu(S_j)^{-1/2} h_j e_{(A_{j}^{-1}x_{j,k})}  : k \in K, j \in J \}$
forms a frame of $L^2(\R^d)$ with lower frame bound $mp\mu(\QQ)^{-1}$
 and
upper frame bound  $MP\mu(\QQ)^{-1}$.
This gives that,
$$
\{ |B_j|^{1/2} h(B_j\omega) e_{x_{j,k}}(B_j\omega) : k \in K, j \in J \}
$$
forms a frame of $L^2(\R^d)$ with frame bounds $mp$ (or $mc$) and $MP$.
The theorem now follows from an application of the inverse Fourier transform.

\end{proof}

\begin {remarks}
\
\begin {enumerate}
\item The set of matrices $\{ A_j \in  GL_d(\R):\,j\in
J \}$ can be
arbitrary and need not have a group structure.
\item The set $\{ A_j \in  GL_d(\R):\,j\in J \}$
can also be chosen to have
a simple structure. For example, $J=\Z^2$, $A_{(i,j)}=R^iD^j$
where $R$ is a rotation and $D$ a dilation matrix, will be used to
construct directional
wavelets. An even simpler example is $J=\Z$, $A_j=A^j$, where
$A$ is an invertible matrix
which gives a construction of wavelet frames on $\R^d$.
\item Note that $h$ does not need to be compactly supported.
\end {enumerate}
\end {remarks}

We will use the
theorem above to construct specific
examples of wavelets, e.g., directional wavelets, isotropic wavelets, etc.

{\bf Interesting particular cases of Theorem~\ref{wavelet}.}
\
\begin{enumerate}
\item $x_{j,k} = x_k\; \forall \; j\in J$.

\noindent
Let $X = \{x_{k} \in \R^d : k \in K\}$ be such that
$\{e_{x_{k}}\chi_{\QQ}, k\in K\}$
is a frame for
$\mathcal{K}_{\QQ}$ with  frame bounds $m$ and
$M$. Then
\begin{equation*}
\{|\det A_j|^{1/2}\psi(A_jx-x_{k}):\, j\in J, k\in K\}
\end{equation*}
forms a wavelet frame of $L^2(\R^d)$ with bounds $mp$ and $MP$. 
\item $A_{j} = A^j \;\forall \; j\in J$, with $A \in GL_d(\R)$.

\noindent
Each of the following sets are wavelet frames of $L^2(\R^d)$ with
bounds
$mp$ and
$MP$:
\begin{eqnarray*}
      & & \{|\det A|^{j/2}\psi(A^jx-x_{j,k}):\, j\in J, k\in K\}, \; \text{and}
\\
      & & \{|\det A|^{j/2}\psi(A^jx-x_{k}):\, j\in J, k\in K\}.
\end{eqnarray*}
\end{enumerate}

\begin {remarks}
\
\begin {enumerate}
\item If the set $\Ss$ is a tiling of $\R^d$ then the
wavelets constructed above are
Shannon-like wavelets, thus not well localized in space. To obtain
well localized space-frequency
wavelets, $\mathcal H$ must be constructed  to
be a smooth partition of unity, e.g., at least $C^1(\R^d)$ as
demonstrated in the examples in Section \ref
{EW}, below.
\item Reconstruction formulas for such wavelet frames are developed in
Section \ref {RF}.
\end {enumerate}
\end {remarks}

\section {Examples of wavelet frames on irregular lattices and with
arbitrary set of dilation
matrices and other transformations}
\label{EW}
\subsection{Density and Separation}
\
To be able to use Theorem \ref {wavelet}
to construct concrete examples of wavelet frames on
irregular grids, we first need to
construct exponential frames $\{e_{x_{k}}\chi_{D}\}$ (also called
Fourier frames) for
$\mathcal {K}_D$. Exponential frames play a central role in sampling
theory for Paley-Wiener
spaces (also known as spaces of band-limited functions).

The density
of a set $X=\{x_k\subset
\R^d: \, k \in K\}$ and separateness of the points in $X$ play a
fundamental role for
finding exponential frames $\{e_{x_{k}}\chi_{D}\}$ for $\mathcal {K}_D$.
\begin {definition}
A sequence $X=\{x_k: k \in K\}$ is separated if
$$
\inf_{k\ne l}\|x_k-x_l\|>0.
$$
\end {definition}
There are many notions for the density of a set $X$. We start with three
definitions that
are due to Beurling.
\begin {definition}
\
\begin {enumerate}
\item A lower uniform density $D^-(X)$ of a separated sequence
$X\subset \R^d$ is defined as
$$
D^-(X)=\lim\limits_{r\to \infty} \frac {\nu^-(r)} {(2r)^d}
$$
where $\nu^-(r):=\min\limits_{y \in \R^d} \card\left({X\cap
(y+[-r,r]^d}\right)$, where $\card(Z)$ denotes the
cardinal of the set $Z$.
\item An upper uniform density $D^+(X)$ of a separated sequence $X$ is
defined as
$$
D^+(X)=\lim\limits_{r\to \infty} \frac {\nu^+(r)} {(2r)^d}
$$
where $\nu^+(r):=\max\limits_{y \in \R} \#\left( {X\cap (y+[-r,r]^d}\right)$.
\item If $D^-(X)=D^+(X)=D(X)$, then $X$ is said to have uniform
Beurling density $D(X)$.
\end {enumerate}
\end {definition}
\begin {remark}
The limits in the definitions of $D^-(X)$ and $D^+(X)$ exist
(see \cite {BW99}).
\end {remark}
As an example, let $X\subset \R$ be separated and assume that there
exists $L>0$ such that  $ |x_k-\frac {k} {d}|\le L$, for all $k \in \Z$.
Then $D^-(X)=D^+(X)=d$.
\
For the one dimensional case, Beurling proved
the following Theorem.
\begin {theorem} (Beurling)
\label {B1}
Let $X \subset \R$ be separated, $a  > 0$ and $\Omega=[-\frac{a} {2},\frac{a}
{2}]$. If $a<D^-(X)$  then
$\{e_{x_{k}}\chi_{\Omega}\}$ is a frame for $\mathcal {K}_\Omega$.
\end {theorem}
This previous result however is only valid in one dimension. For higher
dimensions, Beurling introduced the following notion:
\begin {definition}
The gap $\rho$ of the set $X=\{x_k:k \in K\}$ is defined as
$$
\rho=\rho(X)=\inf \left\{ {r>0: \, \bigcup_{k\in K} B_r(x_k)=\R^d}\right\}
$$
\end {definition}
Equivalently, the gap $\rho$ can be defined as
$$\rho=\rho(X)=\sup_{x\in \R^d} {\inf_{x_k\in
X}|x-x_k|}.
$$

It is not difficult to show that if $X$ has gap $\rho$, then
$D^-(X)\ge \frac {1}{2\rho}$.
For a separated set $X$, and for the case where
$\Omega$ is the ball $B_{r}(0)$ of radius $r$ centered at the origin,
Beurling \cite{Beu66} proved the following result:
\begin {theorem}[Beurling]
\label {beurlingthm}
Let $X\subset \R^d$ be separated, and $\Omega=B_{r}(0)$. If
$r\rho<1/4$, then $\{e_{x_{k}}\chi_{\Omega}\}$
is a frame for $\mathcal {K}_\Omega$.
\end {theorem}
For a very clear exposition of some of the Beurling density results see
\cite{BW99}.

\subsection {Wavelet frames in 1-D}
\label {WF1D}
The construction in the following theorem (which is a particular case
of Theorem~\ref{wavelet}), generalizes a similar result
of \cite{DGM86} to the irregular case. See also \cite{Gro93}.
\begin {theorem}
\label {1dwvlt}
Let $\QQ=[-1,-1/2]\cup [1/2,1]$, $0\le \epsilon < 1/2$, and let $\hat
h_{+}$ be a real valued
function such that $\QQ^1_\epsilon:=\supp \hat h_+\subset
[\frac{1}{2}-\epsilon, 1+\epsilon]$,
$|\hat h_+| \le 1$, and $0<c\le |\hat h_+|$ on $[1/2,1]$. Assume that
for each $j \in \Z$, the sequence
$X_j=\{x_{j,k}\}_{k\in \Z}$ is separated and that
$D^-(X_j)>2^{j+1}(1+\epsilon)$. Then for each $j$,
the set $\{ e_{2^jx_{j,k}}\chi_{\QQ_\epsilon}: \, x_{j,k}\in X_j:
k \in \Z\}$
is a frame of $\mathcal {K}_{\QQ_\epsilon}$, where
$\QQ_\epsilon:=\QQ^1_\epsilon\cup(-\QQ^1_\epsilon)$. If  furthermore
the sets $X_j=\{x_{j,k}\}_{k \in
\Z}$ are chosen such that the frame bounds $m_j$ and $M_j$ satisfy
$\inf_{j} m_j=m>0$ and
$\sup_j M_j=M < +\infty$, then the set
$\{2^{j/2}\psi(2^{j}(\cdot-x_{j,k})):\, j\in \Z, k\in \Z\}$
where $\psi(x)=2Re (
h_+(x))$ is a wavelet frame for $L^2(\R)$.
\end {theorem}
\begin {remark}
The wavelet frame constructed in the theorem above is of the form \\
$\{2^{j/2}\psi(2^{j}(\cdot-x_{j,k})):\, j\in \Z, k\in \Z\}$ which is
slightly different form than the one
constructed in Theorem~\ref{wavelet}. This discrepancy is due to a
convenient choice of the irregular set $X_{j}*= \{x_{j,k}\}$ that we have
adopted in the statement of the theorem above.
\end {remark}
Note that the wavelets constructed in the theorem above  are real and
symmetric. Actually, if one wants $\psi$ with good decay,
$\hat h_+$ can be easily constructed to be $C^r$, $r \geq 1$, even
$C^{\infty}$.

As a corollary, if we choose the sampling sets $X_j$ to be
nested, i.e., $X_j \subset X_{j+1}$, we get
\begin {corollary}
\label {1dwvltc1}
Let $\QQ$, $\epsilon$, and $\hat h_+$ be as in Theorem \ref {1dwvlt}.
Assume that the sequences
$X_j=\{x_{j,k}\}_{k \in \Z}$ are  separated and such that
$x_{j,k}=x_{j+1, 2k}$. If $D^-(X_0)>2(1+\epsilon)$, then for each $j$,
the set $\{e_{2^j x_{j,k}}\chi_{\QQ_\epsilon}: \, x_{j,k}\in X_j\}$
is a frame for
$\mathcal K_{\QQ_\epsilon}$. If  furthermore the frame bounds $m_j$ and
$M_j$ satisfy $\inf_{j} m_j=m>0$ and
$\sup_j M_j=M < \infty$, then the set
$\{2^{j/2}\psi(2^{j}(\cdot-x_{j,k})):\, j\in \Z, k \in \Z\}$ where
$\psi(x)=2Re (h_+(x))$ is a wavelet frame for $L^2(\R)$.
\end {corollary}
\begin {proof}
Since $x_{j,k}=x_{j+1,2k }$ we have that $\#(X_{j+1}\cap [-r,r])\ge 2
\#(X_{j}\cap
[-r,r]-1) $. Thus $D^-(X_{j+1})\ge 2D^-(X_{j})$. But $D^-(X_0)>2(1+\epsilon)$,
therefore $D^-(X_j)>2^{j+1} (1+\epsilon)$. The
corollary then follows directly from Theorem \ref {1dwvlt}.
\end {proof}
     From Theorems \ref {1dwvlt} and \ref {wavelet},
we immediately get the following Corollary.
\begin {corollary}
\label {1dwvltc2}
Let $\QQ$, $\epsilon$, and $h_+$ be as in Theorem \ref {1dwvlt}.
Assume that the set $X=\{x_{k}\}_{k \in \Z}$
is separated and that  $D^-(X)>2(1+\epsilon)$. Then, the set
$\{e_{x_{k}}\chi_{\QQ_\epsilon}: \, x_{k}\in X\}$ is a frame of
$\mathcal {K}_{\QQ_\epsilon}$, and the set of functions
$\{2^{j/2}\psi(2^{j}\cdot-x_k):\, j\in \Z, k \in \Z\}$ where
$\psi(x)=2Re (h_+(x))$ is
a wavelet frame for $L^2(\R)$.
\end {corollary}
\subsubsection {Examples}
\begin {enumerate}
\item {\em Shannon-type wavelet frames }: We use Corollary \ref
{1dwvltc1}, with $\epsilon=0$ and $\hat h_+=\chi_{[1/2,1]}$ to get
wavelet frames of the form
$$\psi_{x_{j,k},j}=\{2^{-j/2}\cos(2^{-j-1} 3\pi  (x-x_{j,k}))\,  \hbox{sinc}
(2^{-j-1}\pi (x-x_{j,k}))\}.$$
These wavelets are not well localized since the decay at $\infty$ is
$O(|x|^{-1})$.
\item {\em Shannon-type wavelet bases }: If we choose $X_j$ such that
$|x_{j,k}-2^jk|\le L_j<\frac {2^j} {4}, \quad \forall\; k \in \Z
$, then by Kadec's $1/4$-Theorem, we immediately get that
$\psi_{x_{j,k},j}=\{2^{-j/2}\cos(2^{-j-1} 3\pi
(x-x_{j,k}))\, \hbox{sinc} (2^{-j-1}\pi (x-x_{j,k})):\, j\in \Z,\ k \in
\Z\}$ constructed above
form a wavelet Riesz basis for $L^2(\R)$.
\item {\em Well localized wavelet frames}: For faster decay of the
wavelet frames, we choose $\hat h_+$ to be a smoother function. Let
$\beta_n=\chi_{[0,1]}\ast\cdots\ast\chi_{[0,1]}$ be the B-spline of
degree $n$ (note that $\supp \beta_n = [0,n+1]$). Let
$\epsilon=1/4$, and  $\hat h_+(\xi)=n\beta_{n-1}\left(%
({\xi-1/4}) n)\right)$. Then we get
a wavelet frame of the form
$$\psi_{x_{j,k},j}=\{2^{-j/2}\cos(2^{-j-1} 3\pi ( x-x_{j,k}))
\,\hbox{sinc}^{n}
(2^{-j}\frac{\pi}{n} (x-x_{j,k}))\}.$$
For this case the wavelets decay as $O(|x|^{-n})$.
\end {enumerate}
\subsection {Examples of wavelet frames in $\R^d$}
\label{examples}
\begin {enumerate}
\item For $\R^2$, let $X=\{x_k: k \in \Z\}$, and
let $Y=\{ y_l: l \in \Z \}$. If $D^-(X) >2$ and  $D^-(Y)>2$, then
using Proposition~\ref {rdframes} below for product frames, the
set
$\{e_{(x_k,y_l)}: (k,l)\in \Z^2\}$ form a frame for $\mathcal
{K}_{[-1,1]^2}$. Let $\QQ:=\{(x,y)\in
\R^2: 1/2\le x^2+y^2\le 1\}$, and
$A=2 I$, then $\R^2=\cup_jA^j\QQ$. We can then use Theorem \ref
{wavelet} to construct wavelet frames
for
$L^2(\R^2)$:
\begin {itemize}
\item {\em Shannon-type radial wavelets}: Let $h=\chi_\QQ$, then $h$
is radial. Thus the function $\psi$ defined
as $\hat \psi= h$ satisfies
$\psi(x,y)=g(r)$, where $r=(x^2+y^2)^{1/2}$. We
then construct the Shannon-like wavelet frame for $L^2(\R^2)$, as in
Theorem \ref {wavelet}. A related construction of  non-separable
radial Shannon-type frame wavelets and
multiwavelets can be found in \cite {PGKKH1}, and \cite {PGKKH2}.
\item {\em Well localized radial wavelets}: To construct wavelet frames with
polynomial decay in space, we let
$h(\xi_1,\xi_2)=n\beta_{n-1}\left({ (\xi^2_1+\xi_2^2-1/4)n}\right)$, and
construct the wavelet frames using Theorem \ref {wavelet} (see Figure
\ref {fig2}).
%%%%%%%%%%%%%%%%%%
\vspace{.2in}
\begin{figure}
\centerline {
\includegraphics[width=50mm]{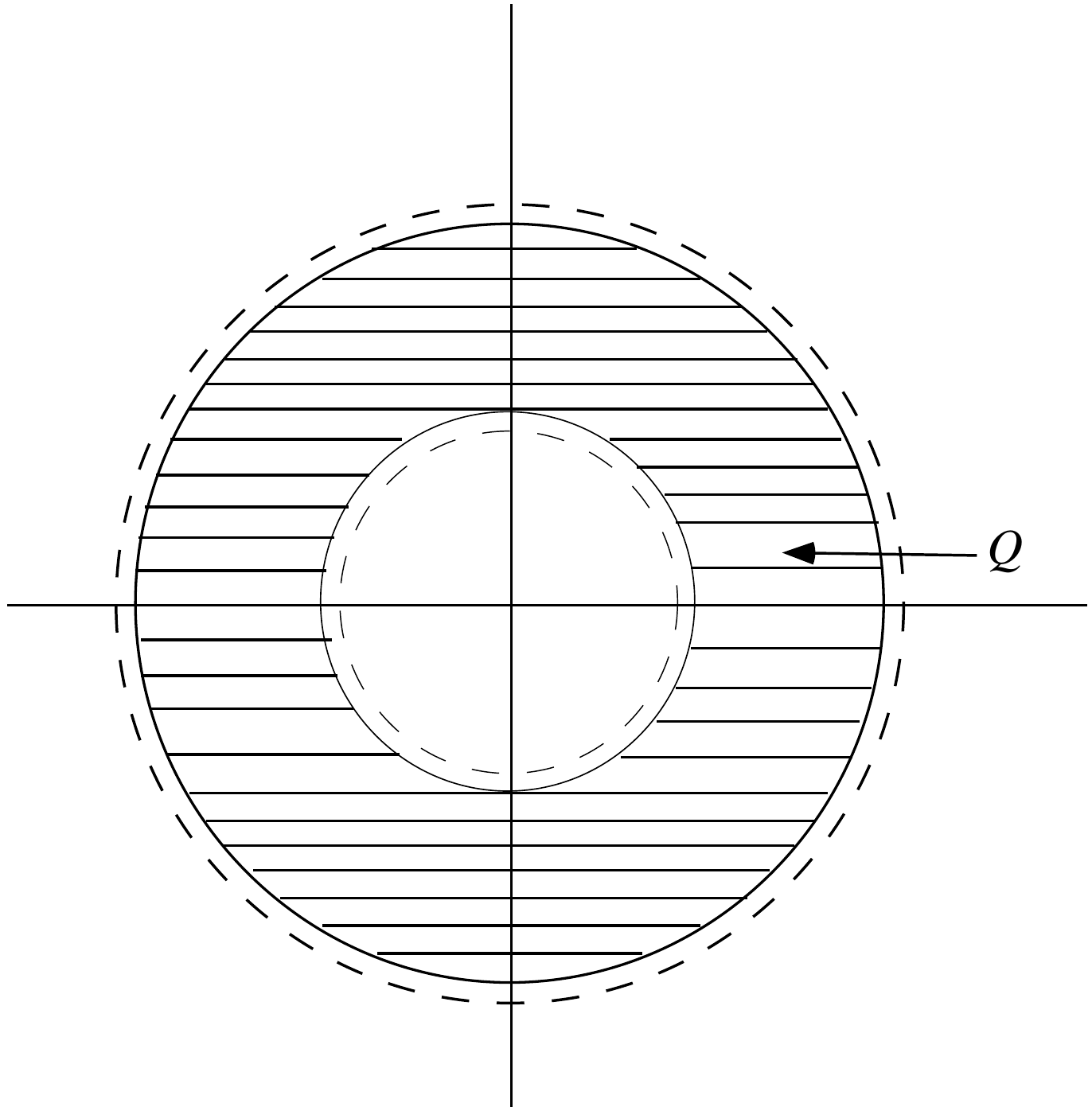}
}
\caption{\label{fig2} Radial wavelet frames that are well localized in space. }
\end{figure}
\vspace{.2in}
%%%%%%%%%%%%%%%%%%
%\ForceWidth{50mm}
%\begin{figure}
%\begin {center}
%\epsfig{file=Fig.2.epsf}
%\epsfig{file=Fig2.ps,width=50mm}
%\caption{\label{fig2} Radial wavelet frames that are well localized in space. }
%\end {center}
%\end{figure}
%%%%%%%%%%%%%%%%%
\end {itemize}
\item The points $\{ (x_k,y_l): k,l \in \Z^2 \}$  lie on an irregular
grid of the form
$X\times Y$. However, we may be interested in points
$Z=\{(x_k,y_l)\subset \R^2: (k,l)\in \Z^2 \}$ that
do not lie on irregular grids of the form $X\times Y$. For this case,
the same constructions above can be used to
form  wavelets frame for $L^2(\R^2)$, as long as the gap
$\rho(X)< \frac {1} {4(1+\epsilon)}$.
\item As in the 1-D examples above, we can also use Corollary \ref
{1dwvltc1} to construct wavelets
on irregular grids satisfying $X_j\subset X_{j+1}$.
\item {\em Directional wavelet frames}: \label{examples-4}
We can easily construct directional wavelet frames as follows:
Let $\QQ_1$ be a region defined
by $\QQ_1=\{(x,y)\in \R^2: x=r\cos(\theta),y=r\sin(\theta), 1/2\le r
\le 1, |\theta|\le \frac
{\pi}{8}\}$, and define $\QQ=(-\QQ_1)\cup\QQ_1$. Let $A=2I$, and $R$
be the matrix of a rotation by an
angle $\pi/4$. Let $\psi$ be such that $\hat \psi=\chi_\QQ$, then we
obtain the wavelet frame for
$L^2(\R^2)$ of the form
$\{\psi_{j,k}=4^{-j_1/2}\psi_j(2^{-j_1}R^{j_2}\cdot -x_{j,k}):
j=(j_1,j_2)\in \Z\times \{0,1,2,3\}, k \in \Z\}$. The index $j_1$ codes for the
resolution of the wavelet, while the
index $j_2$ codes for four possible directions. Thus the wavelet
frame coefficients encode time
scale as well as directional information.  Clearly one can choose any
number of directions and adapt the previous construction. An obvious
modification as shown in
Figure \ref {fig1}, yields wavelet frames with polynomial decay.
%%%%%%%%%%%%%%%%%%
%\vspace{.2in}
\begin{figure}
\centerline {
\includegraphics[width=50mm]{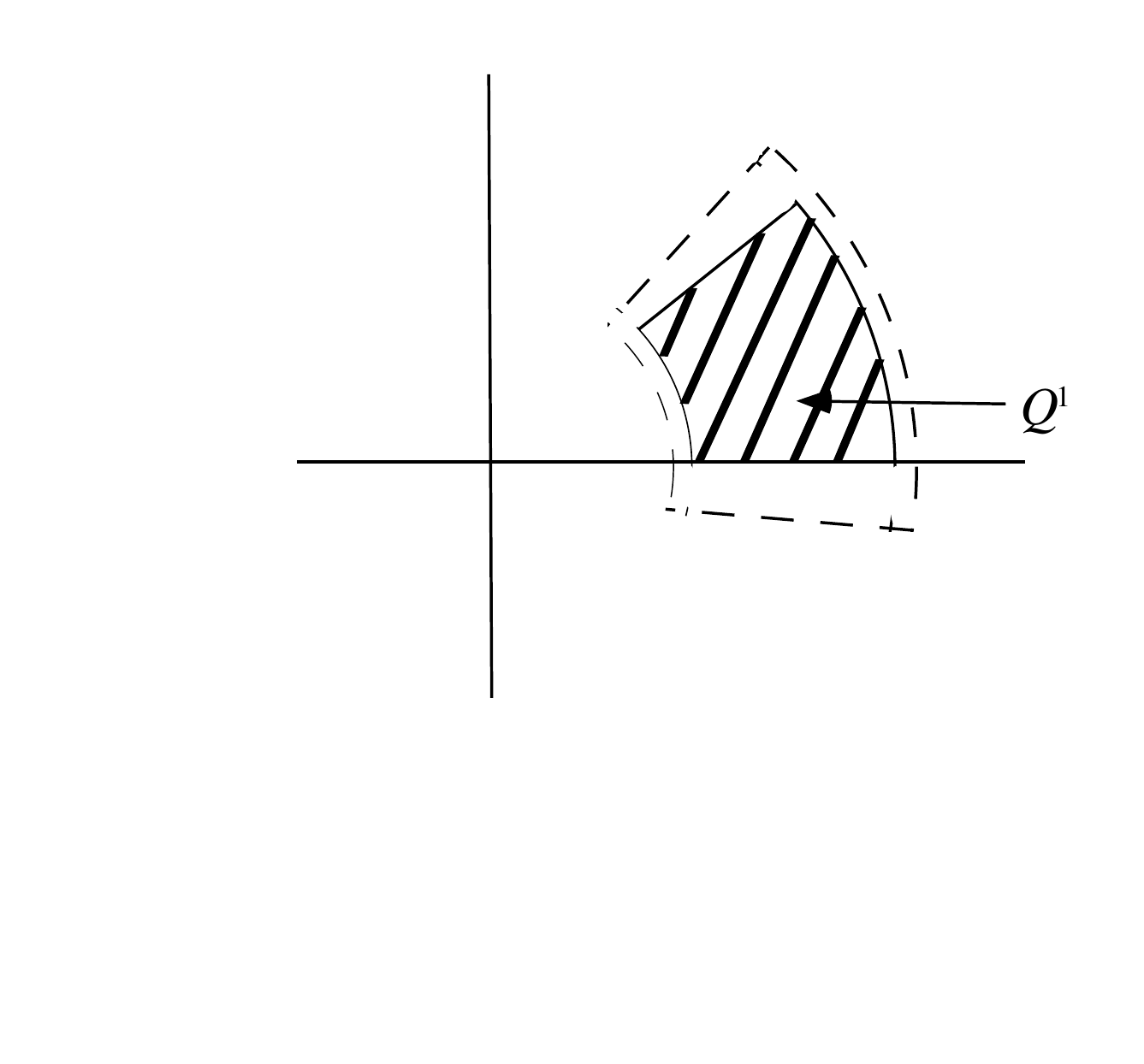}
}
\caption{\label{fig1} Well-localized directional wavelet: The regions
$\QQ=\QQ^1\cup (-\QQ^1)$ and $\QQ_\epsilon=\QQ_\epsilon^1\cup
(-\QQ_\epsilon^1)$ that can be used
to construct well localized directional wavelets. }
\end{figure}
%\vspace{.2in}
%%%%%%%%%%%%%%%%%
Very
nice constructions of smooth
directional wavelet frames on regular grids were obtained before in
\cite{AHNV01, ADHNV03}.
\item {\em Spiral}
In this example we will define a dilation covering by  spiral annulus sectors.

Let $a,b > 1$, and $\Gamma$ the spiral curve defined by
\begin{equation*}
\Gamma(t) = (a^t \cos (bt), a^t \sin (bt)) \quad t \in \R.
\end{equation*}
For $\alpha \in \R$ define $R_{\alpha}$ to be the rotation
of angle $\alpha$ : $R_{\alpha} =
\left[ \begin{smallmatrix} \cos \alpha &- \sin \alpha\\\sin \alpha &
\cos \alpha \end{smallmatrix}\right]$. The curve $\Gamma$ satisfies:
\begin{equation*}
\Gamma(t+\alpha) = a^{\alpha} R_{b\alpha} \Gamma(t).
\end{equation*}
Note that for positive $\alpha$ the matrix $A = a^{\alpha} R_{b\alpha}$
is expansive.

Now we are ready to define the covering elements.
Set $b = 2\pi$ and $\alpha = \frac{1}{m}$, for some integer $m \geq 2$ so
that $A^m = a I_d$. Define the spiral annulus sector $Q = \{x \in \R^2:
x  = \lambda \Gamma(\beta), 1 \leq \lambda \leq a,
0 \leq \beta \leq \alpha\}$ (see Figure \ref {fig3}). So $Q$ is
compact and $\{A^jQ : j
\in
\Z\}$ is a disjoint  covering of $\R^d \setminus \{0\}$.

Choose $\varepsilon > 0$ sufficiently small and $h$ a smooth function that does
not vanish in $Q$  and with support in $Q_{\varepsilon}$. Define
$\hat \psi = h$.
Select a separated set $X=\{x_k\}_{k \in \Z} \subset \R^2$
such that $\rho(X) < \frac {1} {2\text{diam}(Q_{\varepsilon})}$.

The set $\{a^{j/m}\psi(a^{j/m} R_{-2\pi j/m}(x - x_k), k \in \Z, j \in \Z\}$
form a wavelet frame of $L^2(\R^2)$ generated by a single wavelet $\psi$ that
is band-limited, with good decay and directional in frequency.
%%%%%%%%%%%%%%%%%%
\vspace{.2in}
\begin{figure}
\begin {center}
\includegraphics[width=50mm]{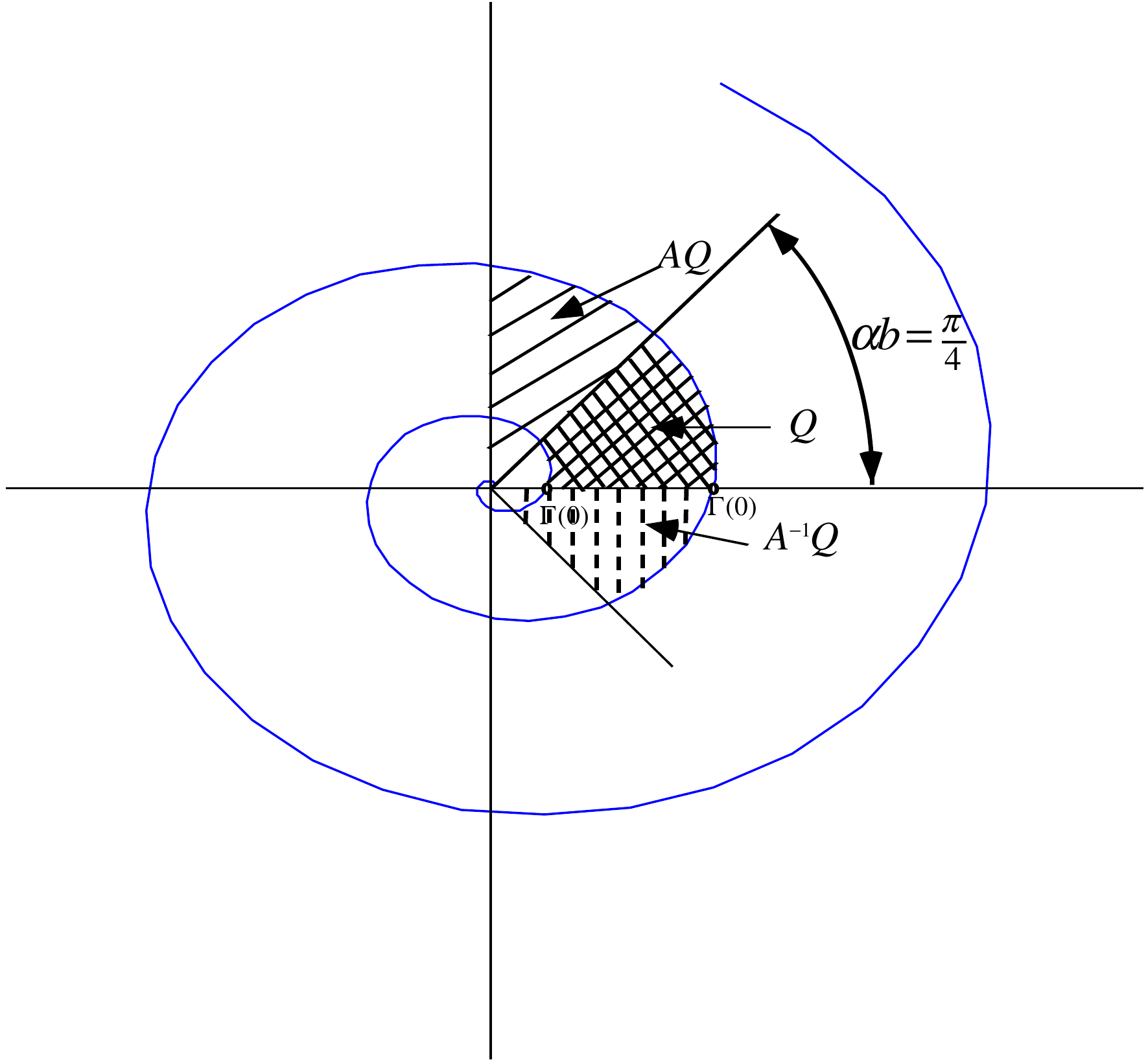}
\caption{\label{fig3} Spiral wavelet frames. }
\end {center}
\end{figure}
%%%%%%%%%%%%%%%%%

\item Obviously, all the constructions above can be generalized to
$\R^d$ for any dimension $d>2$.

\end {enumerate}
\begin {remark}
Some of the wavelet frames may be associated with MRAs. For example,
the so called Shannon wavelet
frame constructed above is associated with the Shannon MRA $V_j=\{ f
\in L^2(\R):
\supp\hat f\subset [-2^{-j-1},2^{-j-1}]\} $, $j\in \Z$. In general
however, the precise relation needs further investigation.
\end {remark}

\section{General results on the construction of time-frequency atoms}
\label{gen-res}
In this section we will develop a general method for constructing
time-frequency frame atoms in several
variables. This construction allows us to construct the previously
introduced wavelet frames on
irregular grids and with arbitrary
dilation matrices or other types of transformations. It also allows
us to construct non-harmonic Gabor
frames on non-uniform coverings of $\R^d$ as
described in Section \ref {CTFA} below.

Let $U$ and $V$ be non-empty open subsets of $\R^d$, $T:U\rightarrow V$
is an invertible $C^1$ map, with $C^1$ inverse $S:V\rightarrow U$,
i.e., $T$ is a $C^1$
homeomorphism with $C^1$ inverse $S:V\rightarrow U$. Define $\alpha
:= \inf_{y\in U} |\det
T^{\prime}(y)|^{-1}$ and
$\beta := \sup_{y\in U} |\det T^{\prime}(y)|^{-1}$, where $T^{\prime}(y)$
denotes the derivative of $T$ at $y$, $T^{\prime}(y) = \left[
\frac{\partial T_i}{\partial y_j} (y)\right]$. We have the
following Proposition :
\begin{proposition}
\label{chvar}
With the above notation, assume that $\alpha$ is positive and $\beta$
is finite, then if
$\{g_j\}_{j \in J}$ is a frame for $\mathcal{K}_V$ with
frame bounds $m$ and $M$, then $\{g_j\circ T\}_{j \in J}$ is
a frame for $\mathcal{K}_U$ with frame bounds $\alpha m$ and $\beta M$.
\end{proposition}
\begin {proof}
Assume that $f \in \mathcal{K}_U$. Then for $j \in J$, by an
application of the change of variables formula
\begin{equation}\label{eq1}
<f,g_{j}\circ T> \ = \int_U f \,\overline{g_{j}\circ T}\
= \int_V (f\circ S)\,|\det S^{\prime}|\, \overline{g_{j}}\ =\
< (f\circ S)|\det S^{\prime}| , g_j>.
\end{equation}
Since $|\det S^{\prime}|$ is finite and bounded away from zero,
$(f\circ S) |\det S^{\prime}|$ is in $\mathcal{K}_V$.

Using that $\{g_{j}\}_{j \in J}$ is a frame of $\KV$
with bounds $m,M$ we have
\begin{equation}
\sum_{j \in J} |<(f\circ S)|\det S^{\prime}|,g_{j}>|^2
\leq M\|(f\circ S)\,|\det S^{\prime}|\,\|^2.
\label{eq2}
\end{equation}
Now, applying again the change of variables theorem, we obtain
\begin{equation*}
\|(f\circ S)|\det S^{\prime}|\,\|^2 = \int_V
|f\circ S|^2 |\det S^{\prime}|^2 = \int_U
|f|^2 |\det (S^{\prime}\circ T)|^2 |\det T^{\prime}|.
\end{equation*}
Since $S = T^{-1}$ we have that $\det[(S^{\prime}\circ T)(y)]
= \left(\det[T^{\prime}(y)]\right)^{-1}$ for all $y \in U$.
Thus
\begin{equation}
\|(f\circ S)|\det S^{\prime}|\,\|^2 =
\int_U
|f|^2 |\det T^{\prime}|^{-1} \leq M\|f\|^2.
\label{eq3}
\end{equation}
That is, from \eqref{eq1}, \eqref{eq2} and \eqref{eq3} we have
\begin{equation*}
\sum_{j \in J} |<f,g_{j}\circ T>|^2 \leq \beta M\|f\|^2.
\end{equation*}
This proves the upper inequality for $\{g_{j}\circ T\}_{j
\in J}$. The lower inequality is obtained in the same way, with the
obvious modifications.
\end{proof}
Proposition \ref{chvar} establishes that $T$
defines an isomorphism
$\Pi_T:\mathcal {K}_U\rightarrow \mathcal {K}_V$ defined by
$\Pi_T(g)=g\circ T$, thereby transforming frames into frames.
By taking $T$ to be a translation or a dilation, we get the well known result:

\begin{corollary}[Translation and dilation of frames]
\label{transl}
Let $\QQ \subseteq \R^d$ be an open subset of $\R^d$.
Let $y \in \R^d$ be any point, and $A \in GL_d(\R)$ an
invertible matrix. We have,
\begin{enumerate}
\item {\bf Translation} \label{coro-1}
$\{g_j\}_{j \in J}$
is a frame for $\KQ$ with bounds $m$ and $M$,
if and only if
        $\{g_{j}(\cdot - y) \}_{j \in J}$ is a frame for $\mathcal{K}_{\QQ+y}$
with the same bounds.
\item {\bf Dilation} \label{coro-2}
$\{g_j\}_{j \in J}$
is a frame for $\KQ$ with bounds $m$ and $M$,
if and only if
        $\{|\det A|^{-1/2} g_{j}(A^{-1} \cdot) \}_{j \in J}$
is a frame for $\mathcal{K}_{A\QQ}$
with the same bounds.
\end{enumerate}
\end{corollary}
\begin {remark}
In fact the corollary remains true if we only assume
that $\QQ$ is measurable.
\end {remark}
\begin{proof}
Part~\ref{coro-1} is a direct application of the Proposition, for the case that
$T(x) = x-y$.

For Part~\ref{coro-2}, and the transformation $T(x) = A^{-1}x$, the Proposition
tells us that $\{g_{j}(A^{-1} \cdot)\}_{j \in J}$ is a frame with
frame bounds
$|\det A|m$ and $|\det A|M$, therefore dividing each function by
$|\det A|^{-1/2}$
we obtain the result.
\end{proof}

For the next theorem, we need to introduce some definitions.
Let $\mathcal{H} = \{h_j:\, j \in J\}$
be a RPU with bounds
$p$ and $P$.  For each $j \in  J$ set
\begin{equation} \label{wj}
W_j = Closure_{L^2}\{\overline h_j f: f \in L^2(\R^d)\}.
\end{equation}

Let $0<c\leq p$, and define
\begin{equation} \label{qj}
Q_j = Q_j(c) = {\{x \in \R^d : |h_j (x)|^2 > c\}}, \quad \text{for
each} \quad j \in J.
\end{equation}
For a given $c$, we discard all those $j$ such that $Q_j$ has measure zero.
Note that if $J_0 = \{j \in J : \mu (Q_j) > 0\}$, then we can only claim that
$\HH_0 = \{h_j\chi_{Q_j}, j\in J_0\}$ is a RPU associated to $ \{\QQ_j\}_{j \in 
J_0}$ with constants
$c$ and $P$.
\begin{theorem}
\label{unionframe}
Let $0< c \le p$.
\begin{enumerate}
\item \label{unionframe-1}
Assume that
$\{g_{j,k}\}_{k\in K}$ is a frame for $W_j$ with
lower and upper frame bounds $m_j$ and $M_j$ respectively.
If $m := \inf_{j} m_j>0$ and
$M := \sup_j M_j < +\infty$, then
$\big\{h_j g_{j,k}: \; j \in J ,  k\in K\big\}$ is a frame
for $L^2(\R^d)$, with frame bounds $pm$ and $PM$.
\item \label{unionframe-2}
Assume  that
for each $j \in J_0$, $\{g_{j,k}\}_{k\in K}$ is a frame for $\KQ_j$ with
lower and upper frame bounds $m_j$ and $M_j$ respectively.
If $m := \inf_{j\in J_0} m_j>0$ and
$M := \sup_{j \in J_0} M_j < +\infty$, then
$\big\{h_j g_{j,k}: \; j \in J_0 ,  k\in K\big\}$ is a frame
for $\mathcal{K}_{\cup Q_j}$, with frame bounds $cm$ and $PM$.
\end{enumerate}
\end {theorem}
\begin {proof}
\
\begin{enumerate}
\item
Given  $f \in L^2(\R^d)$ and denoting by $f_j = \overline h_j f \in W_j$
we will first show that
\begin{equation}\label{fj}
p\|f\|^2  \leq \  \sum_j \|f_j\|^2 \leq P\|f\|^2,
\end{equation}
for
\begin{align*}
p\|f\|^2 & = \ \int p|f|^2 \  \leq   \int \sum_j |h_j|^2 |f|^2  =
\ \sum_j \int |h_j f|^2 = \sum_j \| f_j\|^2,
\end{align*}
where we used dominated convergence for the interchange of the integral
with the sum. The other inequality is analogous.

For each $j \in J$ and $f \in L^2(\R^d)$, we use the fact that $\{g_{j,k}\}_k$
is a frame for $W_j$, and that
$ <f,h_j g_{j,k}> = <\overline h_j f, g_{j,k}>, $
to obtain
\begin{align*}
m \|h_j f\|^2  & \leq
m_j \|h_j f\|^2  \leq
       \ \sum_k |<f,h_j g_{j,k}>|^2 \\
       & \leq
\ M_j \|h_j f\|^2 \ \leq \ M\|h_j f\|^2.
\end{align*}
So summing over $j$
\begin{align*}
pm\|f\|^2 & \leq \ m\sum_j\|f_j\|^2 \ \leq
\ \sum_j\sum_k |<f,h_j g_{j,k}>|^2 \\
& \leq \ M\sum_j\|f_j\|^2 \ \leq \ PM\|f\|^2.
\end{align*}
\item
For the second case we observe that the set
$\lambda = \{\lambda_j = h_j\chi_{Q_j}, j \in J_0\}$ forms a
RPU associated to $\{ Q_j, j \in J_0\}$ with
bounds $c$ and $P$, for if $x \in \cup Q_j$,
\begin{align*}
c & \leq c\, \card (\{j \in J_0: x \in Q_j\}) \leq
\sum_{j \in J_0: x \in Q_j} |h_j(x)|^2 \\
& =
\sum_{j \in J_0} |h_j(x)|^2 \chi_{Q_j}(x) \leq
\sum_{j \in J_0} |h_j(x)|^2 \leq P.
\end{align*}
Furthermore, by hypothesis,
$$ \KQ_j = \{\overline \lambda_j f: f \in L^2(\R^d)\}, \ j \in J_0.
$$
These subspaces $\KQ_j$ correspond to the subspaces $W_j$ defined in \eqref{wj}
for the RPU $\lambda$,
and therefore we can use the previous result (applied to $L^2(\cup Q_j)$
instead of $L^2(\R^d)$) to conclude that
$\{\lambda_j g_{j,k}: j \in J, k \in K\}$ is a frame for
$\KK_{\cup Q_j} = L^2(\cup Q_j)$.
The proof is complete by noting that
$\lambda_j g_{j,k} = h_j g_{j,k}, j \in J_0$.
\end{enumerate}
\end{proof}

%%%%%%%%%%%%%%%%%%%%%%%%%%%

\begin{remarks}
\
\begin{itemize}
\item Note that in the previous theorem,
instead of choosing a frame for the subspaces $W_j$, we could
have chosen any collection of functions of $L^2(\R^d)$ that
form an outer frame for $W_j$.
\item If $h$ is a bounded function and $\QQ= \text{Supp}\, h$, then it is easy
to see that $ \text{closure}_{L^2}\{h f: f \in L^2(\R^d)\}= \mathcal{K}_{\QQ}$
if and only if
$\mu(\QQ) = \mu (\{x: |h(x)| > 0\})$.
So  the spaces $W_j$ defined in \eqref{wj} will coincide
in most of the cases with
$\mathcal{K}_{\text{Supp}(h_j)}$.
\end{itemize}
\end{remarks}
As a very important particular case of the previous theorem, we have the
following Corollaries.

\begin{corollary}
\label{unionframe-cor}
Let $\Ss = \{ S_j \subset \R^d: \; j \in J\}$
be a family of subsets of $\R^d$, not necessarily disjoint,
and let $\mathcal{H} = \{h_j\}$ be a RPU with constants
$p,$ and $P$ not necessarily associated to $\Ss$.
Assume that
$\{g_{j,k}\}_{k\in K}$ is a frame for $\mathcal{K}_{S_j}$ with
lower and upper frame bounds $m_j$ and $M_j$ respectively.
If $m := \inf_{j} m_j>0$ and
$M := \sup_j M_j < +\infty$, then
\begin{enumerate}
\item \label{unionframe-cor-1}
If $\HH$ is associated to $\Ss$ (i.e. $\supp h_j \subseteq S_j,
j \in J$), then
$\big\{h_j g_{j,k}: \; j \in J ,  k\in K\big\}$ is a frame
for $\mathcal {K}_{\cup_{S_j}}$, with frame bounds $pm$ and $PM$.
\item \label{unionframe-cor-2}
If instead
$|h_j(x)|^2 > c\ \forall x \in S_j, j\in J$, then also
$\big\{h_j g_{j,k}: \; j \in J ,  k\in K\big\}$ is a frame
for $\mathcal{K}_{\cup S_j}$, with frame bounds $cm$ and $PM$.
\end{enumerate}
\end {corollary}
\begin{proof}
Denote by $P_{W_j}$ the orthogonal projection on the subspace $W_j$
defined in \eqref{wj}.
Since $\{g_{j,k}: k \in K\}$ forms a frame for $\mathcal{K}_{S_j}$
with bounds  $m_j$ and $M_j$, then $\{ P_{W_j}g_{j,k} : k \in K\}$
forms a frame for $W_j$ with the same bounds.
So, by Theorem~\ref{unionframe},
$\{h_j P_{W_j}g_{j,k} :j \in J, k \in K\}$
forms a frame of $\mathcal{K}_{\cup S_j}
$ with frame bounds $mp$ and $MP$.
We obtain part $(1)$ of the Corollary observing that:
$$ <f,h_j P_{W_j}(g_{j,k})> =  < \overline h_j f, P_{W_j}(g_{j,k})>
= <\overline h_j f , g_{j,k}> =  <f, h_j g_{j,k}> .$$

The second part is a consequence of the fact that
    $\{h_j \chi_{S_j} : j \in J \}$
is a RPU associated to $\Ss$ with bounds $c$ and $P$
and that $\QQ_j(c)= S_j$.
Therefore we can apply the second part of Theorem \ref{unionframe}.
\end{proof}
%%%%%%%%%%%%%%%%%%
\begin {remarks}
\
\begin{itemize}
\item For the first case,
the set $\big\{h_j g_{j,k}: \; k \in K\big\}$ is not
necessarily  a frame for
$\mathcal {K}_{S_j}$ or even for $W_j\subset \mathcal {K}_{S_j}$,
even though $\big\{h_j g_{j,k}: \; j \in J, k\in K\big\}$
is a frame for $\mathcal {K}_{\cup_{S_j})}$.
\item In the second case, the
subspace $W_j$ is contained in $\KS_j$.
\item Note that if $\HH$ is a regular partition of unity, then the frame bounds
for the frame constructed in the Corollary, are $m$ and $M$.
\end{itemize}
\end{remarks}

%%%%%%%%%%%%%%%%%%%%%%%%%%%%%%%%%%%%%%%%%%%%%%%%%%
% CHANGES AFTER REFEREE OF DENSITY

\begin{corollary}
\label{new-cor}
Let $\mathcal{H} = \{h_j: \; j \in J\}$ be a RPU with constants
$p,$ and $P$ and $\{ S_j \subset \R^d: \; j \in J\}$ and $\{ Q_j \subset \R^d: \; j \in J\}$be coverings of $\R^d$ such that $Q_j \subset S_j$ and $|h_j(x)|^2 \geq c\ {\rm a.e. }\  x \in Q_j$ for all $j \in J$ and some constant $c > 0$. 
Assume that
$\{g_{j,k}\}_{k\in K}$ is a frame for $\mathcal{K}_{S_j}$ with
lower and upper frame bounds $m_j$ and $M_j$ respectively.
If $m := \inf_{j} m_j>0$ and
$M := \sup_j M_j < +\infty$, then
$\big\{h_j g_{j,k}: \; j \in J ,  k\in K\big\}$ is a frame
for $ L^2(\R^d)$, with frame bounds $cm$ and $PM$.
\end{corollary}
\begin{proof}
As in the proof of the Theorem, we have the inequalities
\begin{equation} \label{fj-1}
p\|f\|^2  \leq \  \sum_j \|\overline{h}_jf\|^2 \leq P\|f\|^2,
\end{equation}
On the other side,
$$ <f,h_j g_{j,k}> =  < \overline h_j f, g_{j,k} \chi_{S_j}>
= <\overline h_j f \chi_{S_j}, g_{j,k}>,$$
and then
$$
\sum_j \sum_k | <f,h_j g_{j,k}>|^2 \leq M \sum_j \|\overline h_j f \chi_{S_j}\|^2_2 
\leq M \sum_j \|\overline h_j f\|^2_2 \leq MP \|f\|^2_2,
$$
and
$$
\sum_j \sum_k | <f,h_j g_{j,k}>|^2 \geq m \sum_j \|\overline h_j f \chi_{S_j}\|^2_2 
\geq m c\sum_j \|f \chi_{Q_j}\|^2_2 \geq mc \|f\|^2_2.
$$

\end{proof}
%%%%%%%%%%%%%%%%%%%%%
The following proposition is a direct application of Fubini's theorem,
and allows us to construct frames in a product space. By induction the
Theorem can be generalized to hold for any finite number of factors.

\begin{proposition}[Product of Frames]
\label {rdframes}
Let $E_1$ and $E_2$ be measurable subsets
of $\R^{d_1}$ and $\R^{d_2}$ respectively, and let
$\{h_{j }\}_{j  \in J}$ and
$\{g_{k}\}_{k \in K}$ be frames for
$\mathcal{K}_{E_1}$ and
$\mathcal{K}_{E_2}$ with frame bounds $m_1,M_1$ and
$m_2,M_2$. Then $\{h_{j}g_{k}\}_{j, k \in
J\times K}$
is a frame
of $\mathcal{K}_{E_1\times E_2}$ with frame bounds
$m=m_1 m_2, M=M_1 M_2$.
\end{proposition}
\begin {proof}
For any function $f$ in $\mathcal{K}_{E_1\times E_2}$,
\begin {eqnarray*}
\lefteqn{\sum_{j \in J} \sum_{k \in K} \left|
{\int_{E_1\times E_2}
        f(x_1,x_2) \overline{h_{j}(x_1)}
        \overline{g_{k}(x_2)} dx_1 dx_2}\right|^2 =} \\
& &\sum_{j \in J} \left( \sum_{k \in K}
\left|{\int_{E_2}\left( {\int_{E_1}
        f(x_1,x_2) \overline{h_{j }(x_1)} dx_1 }\right)
        \overline{g_{k}(x_2)} dx_2}\right|^2 \right) \\
&  &\ge m_2\int_{E_2}\sum_{j  \in J} \left| \int_{E_1}
        f(x_1,x_2) \overline{h_{j }(x_1)} dx_1
        \right|^2 dx_2 \\
&  &\ge m_1 m_2\int_{E_2} \int_{E_1}
        | f(x_1,x_2) |^2  dx_1  dx_2 ,
\end {eqnarray*}
which yields the lower frame bound $m=m_1m_2$. The upper
frame bound $M=M_1M_2$ can be obtained in a similar fashion.
\end {proof}

%%%%%%%%%%%%%%%%%%%%%%%%%%%%%%

\subsection{Construction of time-frequency atoms on arbitrary
irregular grids and with arbitrary
dilation matrices and other transformations}
\label {CTFA}
We now particularize our previous results to frames of the form
$\{he_{x_{j,k}}, j \in J, k \in K\}$, where $h$ is a fixed function.
Using the Fourier transform, these types of frames
allow us  to construct
wavelets $\{|\det A_j|^{1/2}\psi(A_jx-x_{j,k}): j \in J, k \in K\}$ with
translates on the arbitrary
irregular grid
$X$ and with an arbitrary countable
family of invertible $d\times d$
matrices $\{A_j \in GL_d(\R): \; j \in J\}$ (cf. Theorem~\ref{wavelet}).
First as a particular case of  Corollary \ref{unionframe-cor},
we obtain the following Proposition.

\begin {proposition}
\label {genexpdil}
Assume that $\Ss = \{S_j: \; j\in J\}$ forms a covering of $\R^d$, and let
$\HH=\{h_j : \, j\in J\}$ be a RPU
with bounds $p$ and  $P$ associated to $\Ss$.
Assume also that  $\{{\mu (S_j)^{-1/2}}e_{x_{j,k}}\chi_{S_j}: \; k\in K\}$
is a frame for
$\mathcal{K}_{S_j}$ with  lower and upper frame bounds $m_j$ and
$M_j$ respectively.
If $m := \inf_{j} m_j>0$ and
$M := \sup_j M_{j} < +\infty$, then
$$\big\{{\mu (S_j)^{-1/2}}h_j e_{x_{j,k}}: \; {j \in J, k \in
K}\big\}$$
is a frame for $L^2(\R^d)$ with frame bounds $mp$ and $MP$. 
\end {proposition}
\begin {remark}
If $x_{j,k}$, $S_j$, and $h_j$ are chosen such that
$x_{j,k}=\alpha k, \, k \in \Z^d$, $\alpha \in \R$,
$S_j=S+j,
\,
j \in \Gamma$ and
$h_j=h(\cdot +j), \, j \in \Gamma$, where $\Gamma$ is
a lattice in $\R^d$,
then we obtain the standard Gabor or Weyl-Heisenberg frames. Thus in general,
the construction above can be viewed as non-harmonic Gabor frames with
variable windows $h_j$.
\end {remark}
The following wavelet frame
construction is a direct application of the previous proposition.
Taking $\widehat \psi_j=h_j$, the set
$\{\psi_{j,x_{j,k}}=\psi_j(\cdot -x_{j,k}): j\in J, k \in K\}$ is a
wavelet frame for $L^2(\R^d)$ with frame bounds $mp$ and
$MP$.
Therefore, there exists a dual frame $\tilde \psi_{j,x_{j,k}}$ such that
$$
f=\sum\limits_{j\in J}\sum\limits_{x_{j,k}} <f,\tilde
\psi_{j,x_{j,k}}>\psi_{j,x_{j,k}} \quad \forall
\, f \in L^2(\R^d).$$

We can also use Proposition \ref{genexpdil}
and Beurling Theorem \ref {B1},
to obtain a non-harmonic Gabor frame as the following example shows:
\begin {example}
\label {GF}
       [Non-harmonic Gabor frames]
Let $\QQ=[1,3]$ and Let
$\beta_3=\chi_{[0,1]}\ast\chi_{[0,1]}\ast\chi_{[0,1]}\ast\chi_{[0,1]}$
be the B-spline of  degree
$3$. Clearly, $\R= \cup_{j\in\Z}\QQ+j$. Let
$X_j=\{x_{j,k}\}_{k \in \Z}$ be sets in $\R$ chosen such that for each $j$,
$D^-(X_j)>2$ and such that the frame bounds
$m_j$ and $M_j$ satisfy  $\inf_{j} m_j=m>0$ and
$\sup_j M_j=M < +\infty$. Then using Proposition \ref{genexpdil}
and Beurling Theorem \ref {B1},
we obtain a non-harmonic Gabor frame of the form $\{\beta_3(
t-j)e^{-i2\pi x_{j,k}\omega}: j \in \Z, k \in
\Z\}$. Obviously, we can also use a non-uniform partition of $\R$ and
get  generalized non-harmonic
Gabor frames as discussed earlier.
\end {example}

\begin{remark}
If
$\QQ \subset \R^d$  be a
measurable subset of $\R^d$, and if the family of
sets $\Ss$ is defined by means of expanding or contracting $\QQ$, then
we obtain the theorems of Section~\ref{WC}
as particular cases of Proposition~\ref{genexpdil}.
\end{remark}

\subsection {Construction of exponential frames}
By Corollaries~\ref{transl}~and~\ref{unionframe-cor}, we can
construct a frame for
$\mathcal {K}_\QQ \subset \R^d$ starting from a frame for $\mathcal
{K}_D$ where $D$ is {\em
any} subset of $\R^d$ with nonempty interior. Hence to build  a frame
for $\mathcal {K}_\QQ$ it
is enough to start with a frame for $\mathcal {K}_U$, where  $U$ is
an open disk.
Specifically, since any bounded  measurable set can be covered by a
finite number of
translates of $U$, we can use Corollary~\ref{transl}(1) to find a
frame for $\mathcal {K}_\QQ$.
We can also expand $U$ until it covers $\QQ$ and use
Corollary~\ref{transl}(2). This shows
that there are many ways to construct frames for $\mathcal {K}_\QQ$
starting from a frame for
$\mathcal {K}_U$. Obviously, the particular construction will depend
on the application.

We will describe now two particular constructions.
\begin{enumerate}
\item
Given a frame for $K_D$ where $D$ is a measurable
subset of $\R^d$, we will construct a frame for $K_{\QQ}$.
Let $\Gamma$ be a regular lattice in
$\R^d$ (i.e. $\Gamma = R\Z^d$, where
$R$  is an invertible
$d\times d$ matrix with real entries),
and let $D$ be a measurable subset of $\R^d$ such that $\R^d=
\bigcup\limits_{\gamma
\in \Gamma} (D+\gamma)$ with a finite covering index.
Let $\QQ$  be a measurable subset of $\R^d$ and
define $\QQ_{\gamma} := \QQ\cap (D+\gamma)$ for  $\gamma \in \Gamma$.
Let $\Delta :=\{\gamma \in \Gamma:\; \mu(\QQ_{\gamma}) > 0\}$. By
Corollary \ref {transl}, if
$\{e_{x_k}\chi_D: k \in K\}$ is a frame for $\mathcal{K}_D$,  then
$\{e_{x_k}\chi_{D+\gamma}: k \in K\}$ is a frame
for
$\mathcal{K}_{D+\gamma}$. Hence, $\{e_{x_k}\chi_{\QQ_\gamma}: k \in
K\}$ is also
a frame for
$\mathcal{K}_{\QQ_\gamma}$. Therefore, by Corollary \ref {unionframe-cor},
$\big\{e_{x_k}
\chi_{\QQ_{\gamma}}:\, \gamma \in \Delta, k \in K\big\}$
is a frame for $\mathcal{K}_{\QQ}$.
As an example, when $\QQ$ is a measurable subset of $\R^d$, $D =
[0,1]^d$, and $\Gamma = \Z^d$, we
have that if $\{e_{x_k}\chi_{[0,1]^d}: k \in K\}$  is a frame for
$L^2([0,1]^d)$, then
        $\big\{e_{x_k} \chi_{\QQ_{\gamma}}:\, \gamma \in \Delta, k \in
K\big\}$ is a
frame for $\mathcal{K}_{\QQ}$
(recall that $\QQ_{\gamma} = \QQ \cap ([0,1]^d + \gamma)$).

It is easy to see that the covering requirement
$\R^d= \bigcup\limits_{\gamma
\in \Gamma} (D+\gamma)$ in the previous construction is not
restrictive. Specifically, we only need
$\R^d= \bigcup\limits_{\gamma
\in \Gamma} (\alpha D+\gamma)$, where $\alpha$ is any positive real
number. Furthermore, the
construction remains valid if for each $\QQ_{\gamma}$ we choose a
{\em different} set
$\{x_{\gamma,k}\}$ such that $E_{\gamma} := \{e_{x_{\gamma,k}}\chi_D:
\gamma \in \Delta, k \in K\}$
is a frame for $\mathcal{K}_D$,
and $m := \inf_{\gamma} m_{\gamma} > 0$ and
$M := \sup_{\gamma} M_{\gamma} < +\infty$,
where $m_{\gamma}$ and $M_{\gamma}$ are the lower and
upper frame bounds of $E_{\gamma}$. Finally, note that if $\QQ =
\R^d$ we obtain a frame for
$L^2(\R^d)$.
\item If $\QQ$ is bounded we can construct frames for $K_{\QQ}$ using
Theorem~\ref{beurlingthm}. Let $\delta = \text{diam}(\QQ)$, and
$x_0 \in \R^d$ such that $\QQ\subseteq B(x_0,\delta)$. Let $X =
\{x_k, k \in K\}$ be such that $\rho(X) < \frac{1}{4\delta}$.
Then using Beurlings Theorem (\ref{beurlingthm}) we obtain that
$\{e_{x_k} \chi_{B(0,\delta)}, k \in K\}$ is a frame of $K_{B(0,\delta)}$.
So, using Corollary~\ref{transl}(1),
$\{e_{x_k} \chi_{B(x_0,\delta)}, k \in K\}$ is an outer frame of
$K_{\QQ}$, and $\{e_{x_k}
\chi_{\QQ}, k \in K\}$ is a frame for $K_{\QQ}$.

\end{enumerate}

\subsection{Existence and construction of Riesz partitions of unity}

In view of the previous results,
we will be interested in constructing particular kinds of Riesz
partitions of unity associated to special coverings of the
space. The next results provide the necessary tools to accomplish this task.

If $A$ is a $d \times d$ matrix, we will say that $A$ is {\em expansive}, if
$|\lambda| > 1$ for every eigenvalue $\lambda$ of $A$.

We will use the following known result (see
for example \cite{HJ91},pg.~297).

\begin{lemma} Let $B$ be in  $\C^{d \times d}$ and $\varepsilon > 0$.
There exists a matrix norm $\mnorm{\cdot}$ such that
\begin{equation*}
s(B) \leq  \mnorm {B} \leq s(B) + \varepsilon,
\end{equation*}
where $s(B)$ is the spectral radius of the matrix $B$,
and there exists a norm $\|\cdot\|$ in $\C^d$ such that
\begin{equation*}
\|  Bx \| \leq  \mnorm {B}\;\|x\|.
\end{equation*}
\end{lemma}

As a consequence of this lemma, if $A$ is an expansive matrix, then
there exists a norm $\|\cdot \|$ in $\C^d$ such that
\begin{equation*}
\|A^{-1}x\| \leq c \|x\| \quad 0 < c < 1,
\end{equation*}
and therefore
\begin{equation*}
\|Ax\| \geq c^{\prime} \|x\| \quad   c^{\prime} > 1.
\end{equation*}
In particular for every $x \in \C^d$,
\begin{equation}
\label{exp-mat}
\lim_{j \rightarrow \infty} \|A^{-j} x\| = 0 \quad \text{and} \quad
\text{if}\ x \ne 0 \quad
\lim_{j \rightarrow \infty} \|A^{j} x\| = +\infty.
\end{equation}

\begin{proposition}
\label{prop**}
Let $A$ be a $d\times d$ expansive matrix
and $V \subset \R^d$ a bounded set such that
{\renewcommand{\theenumi}{\roman{enumi}}
\begin{enumerate}
\item there exists $\varepsilon > 0$ such that $B(0,\varepsilon)\cap V =
\emptyset$.
\label{prop**1}
\item $\bigcup_{j\in \Z} A^j V = \R^d \setminus \{0\}$.\label{prop**2}
\end{enumerate}}
Then $\rho_{A,V}$, the covering index of the family $\{A^jV\}_{j \in \Z}$,
is finite, i.e.
there exists an integer $n \ge 1$ such that $1 \leq \rho_{A,V} \leq n$.
\end{proposition}

\begin{proof}
For $x \in \R^d$ define
\begin{equation*}
z^{+}_{x} = \{j \geq 0 : x \in A^jV\} \quad \text{and}\quad
z^{-}_{x} = \{j < 0 : x \in A^jV\}.
\end{equation*}
Since $\rho_{A,V} \leq \supess_x(\card(z^{+}_x)) + \supess_x(\card(z^{-}_x))$,
it is enough to prove that
$\card(z^{+}_x)$
and $\card(z^{-}_x)$ are uniformly bounded in $\R^d$.
We will see that $\card(z^{+}_x)$ is uniformly bounded. A similar
argument proves
the  claim for $\card(z^{-}_x)$.

It is easy to see that
\begin{equation*}
\card(z^{+}_x) = \card(z^{+}_{A^sx}) \quad \forall\ s \in \Z, \quad \forall \
x \in \R^d.
\end{equation*}
Thus, if  $\card(z^{+}_x)$ is bounded on $V$, then, by (\ref{prop**2})
$\card(z^{+}_x)$ is bounded  in $\R^d \setminus \{0\}$ with the same bound.
Now, using (\ref{prop**1}), we see that there exist $0 < c_1 < c_2 < +\infty$
such that $c_1 \leq \|x\| \leq c_2$, for all $x \in V$. Fix $x \in V$.
If $j \in  z^{+}_x$, then there exists $v \in V$ such that $x = A^jv$ and
\begin{equation*}
c_2 \geq \|x\| = \|A^j v\| \geq \frac{\|v\|}{\|A^{-j}\|} \geq
\frac{c_1}{\|A^{-j}\|},
\end{equation*}
and therefore $\|A^{-j}\| \geq \frac{c_1}{c_2}$. But
since $\|A^{-j}\| \rightarrow 0$, necessarily there exists $j_0$ such
that for every $s \geq j_0, {\|A^{-s}\|} \leq \frac{c_1}{c_2}$.
Hence for every $x \in V$, $\card(z^{+}_x) \leq j_0$.
\end{proof}
The next proposition shows the construction of a large class of Riesz
partitions of
unity, for families of sets obtained by dilation of a compact set.
\begin{proposition}\label{prop2}
Let $\QQ \subset \R^d$ be a compact set and $A$ a $d\times d$ expansive matrix
such that
{\renewcommand{\theenumi}{\roman{enumi}}
\begin{enumerate}%[roman]
\item $0 \not\in \QQ$ \label{prop2-1}
\item \label{prop2-3}$\bigcup_{j \in \Z} A^j \QQ =\R^d \setminus \{0\}$.
\end{enumerate}}
Let $h$ be any measurable function, and $0<c_1 \leq c_2 < +\infty$ some
constants such that
{\renewcommand{\theenumi}{\alph{enumi}}
\begin{enumerate}
\item $0 \leq |h|^2 \leq c_2$
\item \label{prop2-5} $0 < c_1 \leq |h|^2$ on $\QQ$
\item $ h = 0$ on $\R^d \setminus \QQ_{\varepsilon}$, where
$0 < \varepsilon < d(0,\QQ)$ and $\QQ_{\varepsilon} \equiv
\{x \in \R^d : d(x,\QQ) \leq \varepsilon\}.$
\end{enumerate}}
Then the family of functions $\{h_j(\cdot) = h(A^{-j}\cdot)\}_{j\in \Z}$ is
a RPU associated to $\{A^j \QQ_{\varepsilon}\}$.
\end{proposition}
\begin{proof}
If $x \not= 0$, by
(\ref{prop2-3}) there exists $j \in \Z$ and $q \in \QQ$ such that
$A^{j}q = x$, so by (\ref{prop2-5}),
\begin{equation*}
c_1 \leq |h(q)|^2 = |h(A^{-j} x)|^2. \quad \text{Thus} \quad
\sum_{s\in \Z} |h(A^{-s} x)|^2 \geq c_1.
\end{equation*}
Now by Proposition~\ref{prop**}, the covering index
$\rho_{A,\QQ_{\varepsilon}}$ of
the family $\{A^j\QQ_{\varepsilon}\}_{j \in \Z}$ is finite.
Using (\ref{prop2-1}) and (\ref{prop2-3}),
we see that $\supp (h_j) \subset A^j\QQ_{\varepsilon}$ and since $0 \leq
|h_j(x)|^2 \leq c_2 \ \forall\ x$, we obtain that
\begin{equation*}
\sum_{j \in \Z} |h_j(x)|^2 \leq \rho_{A,Q_{\varepsilon}} c_2,
\end{equation*}
which  proves the proposition.
\end{proof}

\begin{remark}
This proposition generalizes easily to the case where we replace
for each $j$, $A^j$ by an invertible matrix $A_j$ in such a way that
$\{ A_j \QQ\}_{j\in \Z}$ is a covering of $\R^d \setminus \{0\}$ with
finite index.
\end{remark}

The next Lemma shows, that the assumption of having a compact set
that covers $\R^d$ by dilations is actually necessary, if one wants to
construct a RPU associated to a covering of
$\R^d$ of the form $\Ss = \{A^j V, j \in \Z\}$, where $V$ is a bounded
open set.

\begin{lemma}
If $A$ is invertible, and $V$ is a bounded open set such that
\begin{enumerate}
\item \label{lem-i}
there exists $\varepsilon > 0$ with $B(0,\varepsilon)\cap V = \emptyset$
\item \label{lem-ii}
$\bigcup_{j\in \Z} A^j V = \R^d \setminus \{0\}$
\end{enumerate}
then there exists a compact set $\QQ \subset V$ such that
$\cup_{j\in \Z} A^j \QQ \supset \R^d\setminus \{0\}$.
\end{lemma}
\begin{proof}
It is clearly enough to find $\QQ$ such that $\cup_j A^j\QQ \supset V$.
Let
\begin{equation*}
Q_n = \left\{ x \in V: d(x,\partial V) \geq \frac{1}{n}\right\}.
\end{equation*}
We will prove (by the contradiction) that for some $n \geq 1$, $Q_n$
covers $V$ by dilations by $A$.

Assume that for each $n \geq 1$, there exists $x_n \in V$ such that
$x_n \not\in \cup_j A^j Q_n$. Since $\overline{V}$ is compact,
and $d(x_n,\partial V) < \frac{1}{n}$, there exists a subsequence
$\{x_{n_k}\}$ and
$x \in \partial V$ such that $x_{n_k} \rightarrow x$. By (\ref{lem-i})
and (\ref{lem-ii}) there exists $j_0$ such that $x \in A^{j_0} V$, and
since $A^{j_0} V$ is open, $x_{n_k} \in A^{j_0} V$ for $ k \geq k_0$.
Set $y = A^{-j_0} x$ and $y_k = A^{-j_0} x_{n_k}$, then
$y, y_k \in V$ for $k \geq k_0$. Choose $\varepsilon$ small enough such that
$B(y,\varepsilon) \subset V$, and so $d(y,\partial V) \geq \varepsilon$.
Let $m_0$ be such that if $m \geq m_0$, then $y_m \in
B(y,\frac{\varepsilon}{2})$.

Then for $z \in B(y,\frac{\varepsilon}{2})$, and $v \in \partial V$ we have
\begin{equation*}
\varepsilon - d(y,z) \leq d(y,v) - d(y,z) \leq d(z,v),
\end{equation*}
and thus
\begin{equation*}
\frac{\varepsilon}{2} \leq d(z,\partial V).
\end{equation*}
So $B(y,\frac{\varepsilon}{2})
\subset Q_n$ for all $n$ such that $\frac{1}{n} < \frac{\varepsilon}{2}$.
This contradicts our assumption that $x_{n_m}=A^{j_0}y_{n_m}\notin
\QQ_{n_m}$, since, if $m \geq m_0$
and
$\frac{1}{n_m} < \frac{\varepsilon}{2}$, we
have $y_m \in B(y,\frac{\varepsilon}{2})
\subset Q_{n_m}$, and so $x_{n_m} \in A^{j_0}Q_{n_m}$.
\end{proof}

We show next that it is not difficult to construct bounded sets
that cover $\R^d$  by dilations.

\begin{lemma}
Let $V \subset \R^d$ be a bounded set such that $0 \in V^{\circ}$, and
$A$ an expansive $d \times d$ matrix.
Set $Q =  AV \setminus V$, then $\{A^j Q, j \in \Z \}$ is a covering
of $\R^d \setminus \{0\}$ with finite covering index. Furthermore, if
$V \subset AV$ then the sets $\{A^j Q\}$ are disjoint.
\end{lemma}
\begin{proof}
Choose $\varepsilon > 0$ such that
$B(0,\varepsilon) \subset V$.  Let $x \in \R^d \setminus \{0\}$
be an arbitrary point. By equation~\eqref{exp-mat}
$\lim_{j \rightarrow \infty} \|A^{-j} x\| = 0$,
for some norm $\|\cdot\|$ in $\C^d$.

So there exists a positive integer $n$ such that $\forall j \geq n$,
      $y := A^{-j}x \in B(0,\varepsilon)$,  i.e., $x \in A^j V, \forall
j \geq n$.

Since $V$ is bounded, and $\|x\| > 0$, there exists $j_0 \in \Z$
such that $x \in A^{j_0+1}V$ and $x \not\in A^{j}V, \forall j \leq
j_0$, i.e., $x \in A^{j_0+1}V \setminus A^{j_0}V = A^{j_0} Q$.

The finiteness of the covering index follows from Proposition~\ref{prop**},
and if $V \subset AV$ the disjointness property follows immediately.
\end{proof}

\begin{remark}
\
\begin{itemize}
\item
If we want to obtain coverings of $\R^d\setminus \{0\}$ with compact sets, we
can choose the set $\overline \QQ$ where $\QQ$ is the set
constructed in the previous Lemma. Clearly the sets
$\{A^j\overline \QQ:\, j \in J\}$
cover $\R^d\setminus \{0\}$, and if $\mu (\partial V) = 0$, then
$\{A^j\overline \QQ:\, j \in J\}$
are almost disjoint, i.e.~$\mu (A^j\overline \QQ \cap A^k\overline \QQ) = 0,
j\not=k$.
\item
Even though for {\em any } set $V$
it is true that $A^j (AV \setminus V) \cap A^{j+1}(AV \setminus V) =
\emptyset $,
       in general it is not true that the family $\{A^j (AV \setminus V)\}$ is
disjoint, as the following example shows:
Let $A = \left[ \begin{smallmatrix}0 &-2\\2&0\end{smallmatrix}\right]$
and $V$ the rectangle $[-3,3]\times[-1,1]$.

\end{itemize}
\end{remark}

\subsection{Recipe for constructing smooth wavelet frames using a
single dilation matrix and an
irregular grid}

We can combine the results of this section with
Theorem~\ref{wavelet}, to obtain the following recipe to construct smooth,
wavelet frames of $\R^d$ associated to a single dilation matrix $A$
and an irregular grid. The
wavelets obtained by this method can be constructed to have
polynomial decay of any degree, as
exemplified in Section \ref {WF1D}.

{\bf Recipe:}
\
\begin{itemize}
\item
Select a bounded set  $V \subset \R^d$ such that $0 \in V^{\circ}$
and $\mu (\partial V)= 0$.
\item
Select a function $h$ of class $C^r, r > 0$ such that $|h| \not= 0$ on
$\QQ = \overline{A^TV\setminus V}$ and $\supp h \subset \QQ_{\varepsilon}$
for some small $\varepsilon > 0$.
\item
Consider a set $X = \{x_{k}\}_{k\in K}
      \subset \R^d$, such that $X$ is separated and
$\rho(X) < \frac{1}{2\delta}$, where
$\delta = \text{diam}(\QQ_{\varepsilon})$.
\end{itemize}
Then the following collection is a frame of $L^2(\R^d)$:
\begin{equation*}
\{|\det A|^{j/2}\psi(A^jx-x_{k}):\, k\in K, j\in \Z)\},
\end{equation*}
where $\psi$ is the inverse Fourier transform of $h$.

%%%%%%%%%%%%%%%%%%%%%%%%%%%%%%%%%%%%

\section{Reconstruction formulas}
\label {RF}
We first note that although the set $\{|\det
A_j|^{1/2}\psi(A_jx-x_{j,k}):\, j\in J, k\in K\}$ in Theorem~\ref
{wavelet} is a wavelet frame for
$L^2(\R^d)$, it is  not in general true that for a fixed $j$ the set
$\{\psi_{j,x_{j,k}}(x)=|\det
A_j|^{1/2}\psi(A_jx-x_{j,k}):\,  k\in K\}$ is a frame, unless the
wavelet frame is of Shannon-type,
hence, not well-localized. Thus, for well-localized wavelets, it
appears that the reconstruction of a
function
$f \in L^2(\R^d)$ from the wavelet coefficients
$\{<f,\psi_{j,x_{j,k}}>: j \in J, k \in K\}$ cannot be obtained in a
stable way by first
reconstructing at each level
$j$ and then obtaining $f$ by summing over all levels
$j$.  However what follows shows that this is still possible:

For a matrix $R$, denote by $D_{R}$
the dilation operator $D_{R}(f)= |\det R|^{1/2} f \circ R$.
Using the notation of Theorem \ref{wavelet} we have that

$$
<f,\psi_{j,x_{j,k}}>=<\hat f,h_j D_{B_j}e_{x_{j,k}}> =
<\overline h_j\hat f,D_{B_j}e_{x_{j,k}}>.
$$
This, together with  the assumptions that
$ \{e_{x_{j,k}}\chi_{\QQ}: k \in K\}$
is a frame for $\mathcal {K}_{\QQ}$ for each $j$,
and that $\supp h_j \subset S_j$,
permits to reconstruct $\hat f_j=\overline h_j\hat f$ using a dual frame
$\{\phi_{j,k}: \,  k \in K\}$
    of $\{e_{x_{j,k}}\chi_{\QQ}: k \in K\}$.
Note that, since $D_{B_j}$ is unitary then
$\{D_{B_j}\phi_{j,k}: \,  k \in K\}$ is a dual frame of
$\{D_{B_j} (e_{x_{j,k}}\chi_{\QQ}): k \in K\}$ for $\mathcal{K}_{S_j}$.

Thus, it is always possible to
reconstruct each $f_j$ in a stable way and then obtain $f$
by summing up over all levels $j$. One drawback is
that the dual frame may not be well-localized, since  each
$\phi_{j,k}$ may be discontinuous at the boundary
of $\QQ$. To treat this problem we simply note that
\[
\overline {h_j}\hat f=\sum_k <\hat f,h_j D_{B_j}e_{x_{j,k}}> D_{B_j}
\phi_{j,k}.
\]
Multiplying both sides by $h_j$ we obtain

\begin {equation}
\label {RECFOR}
|h_j|^2\hat f=\sum_k <\hat f,h_j D_{B_j} e_{x_{j,k}}>  D_{B_j}\theta_{j,k},
\end {equation}

where $ \theta_{j,k}= h \phi_{j,k}$. If we choose $h$
to be in $C^r(\R^d)$,
%with $\text{supp } h \subset \QQ$,
    and if
$\phi_{j,k}$ is in $C^r(\QQ)$ ,
then  $\theta_{j,k}$ will be in $C^r(\R^d)$ and therefore decays
polynomially in space.
Hence the partial sums of the series (\ref {RECFOR}) will have
good convergence properties.
We can then sum equation~\eqref{RECFOR} over $j \in J$ to
obtain $\hat f\sum_j|h_j|^2$,
and then divide by
$\sum_j|h_j|^2$ to obtain $\hat f$.

\par {\em Acknowledgments}: We wish to thank Professors Chris Heil,
Yves Meyer, Gestur \'{O}lafsson,
    and Qiyu Sun for   insightful remarks.

%%%%%%%%%%%%%%%%%%%%%%%%%%%%%%%%%%%%

%\bibliographystyle{amsalpha}
%\bibliography{cyu}
\newcommand{\etalchar}[1]{$^{#1}$}
\providecommand{\bysame}{\leavevmode\hbox to3em{\hrulefill}\thinspace}
\providecommand{\MR}{\relax\ifhmode\unskip\space\fi MR }
\providecommand{\MRhref}[2]{%
        \href{http://www.ams.org/mathscinet-getitem?mr=#1}{#2}
}
\providecommand{\href}[2]{#2}

\end{document}